**Title: Behavioral Switching Loss Modeling of Inverter Modules**

Authors: K. Stoyka, R. A. Pessinatti Ohashi, N. Femia

**Acknowledgement:** The results presented in this paper are developed in the framework of the 16ENG08 MICEV Project. The latter received funding from the EMPIR programme co-financed by the Participating States and from the European Union's Horizon 2020 research and innovation programme.

# Behavioral Switching Loss Modeling of Inverter Modules

Kateryna Stoyka, Ricieri Akihito Pessinatti Ohashi and Nicola Femia
DIEM – University of Salerno, Fisciano (SA), ITALY
e-mail:{*kstoyka* ; *rpessinattiohashi* ; *femia*}@unisa.it

*Abstract* — This paper presents a new behavioral model for switching power loss evaluation in phase-shifted full-bridge inverter Power Modules (PoMs). The proposed model has been identified by means of a Genetic Programming (GP) algorithm combined with a Multi-Objective Optimization (MOO) technique. A large set of loss data, evaluated by means of analytical loss formulas, has been considered for the identification of a compact behavioral model. The GP-MOO approach considers the inverter switching frequency, input voltage, duty-cycle and load resistance as model input variables, and the MOSFET gate driver voltage and resistance as parameters influencing the coefficients values of the identified loss formula. The behavioral model loss predictions confirm their reliability for a wide range of operating conditions.

## I. Introduction

The achievement of increasing power density levels is a major concern in power supply units design [1]. Several integrated hardware solutions have allowed reducing the impact of power devices on supply units, such as Power Modules (PoMs), which usually contain a combination of two or more devices of similar technology to realize functional power sub-systems [2]. In recent years, the demand for semiconductors PoMs has considerably increased due to their benefits over discrete design solutions, such as size, cost, time-to-market, reliability and flexible layout [3]. PoMs are available for diodes, MOSFETs, thyristors and Insulated Gate Bipolar Transistors (IGBTs), and can cover a wide range of input voltages and output currents, with a variety of package options, thus resulting advantageous in manifold applications, such as in wireless power transfer systems [4], but also in power factor correction converters, inductive heating systems and solar inverters.

In literature, several models and methods are available for the estimation of power loss of semiconductor devices. Analytical models [5] are usable only if the devices operating conditions correspond to those ones adopted for their characterization. Numerical models [5] can instead provide more accurate loss prediction, but longer simulation times and possible numerical instabilities. When different factors affect the power loss evaluation, some of which difficult to be modelled in a reliable way (e.g., body-diode reverse recovery effect, capacitances charging and discharging), it can be convenient to adopt behavioral models. This is particularly true for semiconductor PoMs, where many devices are possibly integrated, and only their overall loss model as a function of the main operating parameters is of interest, rather than the loss models of each single device. Indeed, the prediction of power systems efficiency, under wide ranges of operating conditions relevant for their application, is a major concern when PoMs are used. Inaccurate analytical loss models, as well as long time-consuming tuning procedures to correct the model parameters estimations are undesired. On the other side, experimental approaches can be considered for accurate loss measurements [7], but they are time consuming too. As an alternative, behavioral loss models for PoMs can be adopted for an effective design, as the operating conditions (e.g. current, voltage, frequency) are the only known quantities.

This paper presents a compact behavioral loss model for Silicon Carbide (SiC) MOSFETs PoMs. Such model has been identified by means of a Genetic Programming (GP) algorithm, in combination with a Multi-Objective Optimization (MOO) technique. Such identification algorithm was previously used for discovering the power loss models of IGBTs and power inductors [8]-[10]. In this paper, the GP-MOO approach is adopted to identify the switching power loss model of a phase-shifted full-bridge inverter PoM. The model considers the inverter switching frequency, input voltage, duty-cycle and load resistance as model input variables, while MOSFETs gate-driver voltage and resistance are used as parameters influencing the model coefficients. The GP-MOO approach has been applied to a large set of switching power loss data, herein generated by means of analytical loss models of devices available in literature [5], emulating experimental loss data sets.

The paper is organized as follows. In Section II, literature power loss formulas valid for the phase-shifted full-bridge inverter are summarized. Section III illustrates the main elements of the GP-MOO approach. In Section IV, the behavioral switching loss model, based on the GP-MOO approach, is discussed. Its predictions are eventually compared to the analytical formulas results for final validation.

## II. Power Loss Modeling for Inverter Module

Fig. 1 shows the circuit schematic of a phase-shifted full-bridge inverter with resonant load represented by the equivalent impedance $\dot{Z}_T$. The phase-shifted control involves a phase-lag in the control signals of the two full-bridge arms [11]. For the first harmonic approximation, the inverter output voltage and current are related to each other as $\dot{Z}_T = \bar{V}_1/\bar{I}_1 = Z_T\angle\phi$.

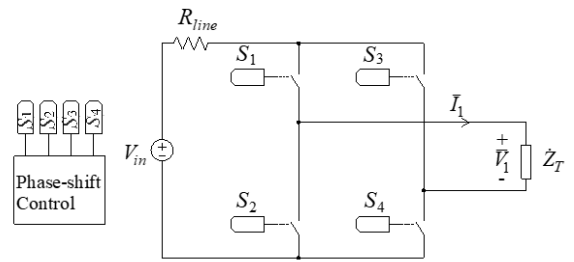

Fig. 1. Phase-shifted full bridge inverter

As a result, square-wave voltage and quasi-sinusoidal current waveforms are achieved at the inverter output, as shown in Fig. 2. Let us consider the half period $[T_s/2, T_s]$, where $\alpha$ is the angle between the zero-crossing current $i_1(t)$ and the falling voltage $v_1(t)$, $\beta$ is the angle between the rising voltage $v_1(t)$ and the zero-crossing current $i_1(t)$, and $D$ is the inverter duty-cycle (fraction of half a period in which the inverter delivers power to the load). It is easily verified that $D = 1-(\alpha+\beta)/\pi$ and $\phi = (\beta-\alpha)/2$ [12]. In perfect resonant matching, the impedance $\dot{Z}_T$ coincides with its real component $R_T$, and the voltage $\overline{V_1}$ and current $\overline{I_1}$ are in phase ($\phi = 0$). This implies that $\alpha = \beta > 0$. Applying Fourier formulas to the voltage $v_1(t)$ provides:

$$V_1 = \frac{4}{\pi} V_{in} \sin\left(\frac{\pi}{2} D\right) \quad (1)$$

Thus, the current magnitude $I_1$ can be obtained from the values of $V_1$ and $R_T$ and used to calculate the switches current values at the transition instants ($I_{1\alpha} = I_1 \sin\alpha = I_1 s_\alpha$; $I_{1\beta} = I_1 \sin\beta = I_1 s_\beta$).

The inverter module total power loss can be evaluated as the sum of conduction and switching losses. The conduction power loss $P_{cond}$ can be evaluated as given in (2):

$$P_{cond} = 2 R_{DS} I_{1,rms}^2 \quad (2)$$

The conduction power loss contribution can be easily estimated from the device drain-source resistance $R_{DS}$ and the inverter rms current $I_{1,rms}$. Conversely, the identification of a behavioral model for the PoM switching power loss term remains the major issue. In fact, the switching power loss $P_{sw}$ is the sum of several contributions, including body-diode loss $P_{bd}$, gate loss $P_{gt}$ and voltage-current overlapping loss $P_{ov}$, where:

$$P_{bd} = 2 f_s t_{dt} V_{SD} I_1 (s_\alpha + s_\beta) \quad (3)$$

$$P_{gt} = 4 f_s V_{dr} Q_g \quad (4)$$

$$P_{ov} = f_s V_{in} \left[ I_{1\alpha} t_{on} + I_{1\beta} t_{off} + C_{oss} V_{in} \right] \quad (5)$$

$$t_{on} = \frac{Q_{gsw} R_{g,on}}{V_{dr} - V_{GS,th} - I_{1\alpha}/g_{fs}} \quad (6)$$

$$t_{off} = \frac{Q_{gsw} R_{g,off}}{V_{GS,th} + I_{1\beta}/g_{fs}} \quad (7)$$

where $Q_g$, $Q_{g,sw}$, $C_{oss}$, $V_{GS,th}$, $g_{fs}$, $V_{SD}$, $t_{dt}$, $V_{dr}$, $R_{g,on}$, $R_{g,off}$ are, respectively, the MOSFETs total and switching gate charge, output capacitance, gate-source threshold voltage, trans-conductance, body-diode voltage, gate signal dead-time, gate-driver voltage, turn-on and turn-off resistances [5][6]. The MOSFET operating temperature influences the values of the parameters $R_{DS}$, $V_{GS,th}$ and $g_{fs}$. Hence, a thermal analysis is required for power loss calculation.

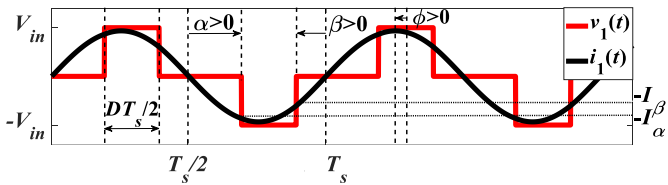

Fig. 2. Voltage and current output waveforms of the inverter module.

Let us consider four SiC MOSFETs in the inverter of Fig. 1, operating at a uniform temperature $T$. Given the ambient temperature $T_a = 25°C$, the $R_{DS}$, $V_{GS,th}$ and $g_{fs}$ values can be updated by means of (8)-(10), and MOSFETs total power loss evaluated accordingly:

$$R_{DS}(T) = R_{DS@25°C}\left[1 + \rho_T(T - T_a)\right] \quad (8)$$

$$V_{GS,th}(T) = V_{GS,th@25°C}\left[1 + \upsilon_T(T - T_a)\right] \quad (9)$$

$$g_{fs}(T) = g_{fs@25°C}\left[1 + \gamma_T(T - T_a)\right] \quad (10)$$

where $\rho_T$, $\upsilon_T$, and $\gamma_T$ are linear thermal coefficients of $R_{DS}$, $V_{GS,th}$ and $g_{fs}$, respectively. From the total power loss $P_{tot} = P_{cond} + P_{sw}$ and the module equivalent thermal resistance $R_{th}$, the new module temperature can be estimated as $T_{new} = T_a + R_{th} P_{tot}$. This thermal routine stops if the normalized difference between the new temperature value and the previous one is lower than a certain threshold (e.g., 1e-4).

In this paper, the thermal-based analytical model presented through formulas (2)-(7) has been used to simulate the inverter switching power loss of a SiC MOSFETs PoM under a wide range of different operating conditions. In particular, given these simulated data, a new compact behavioral switching power loss model has been identified and presented as function of the operating conditions imposed by the application of interest. Two main assumptions have been done. First, the total PoM switching power loss depends on the inverter switching frequency $f_s$, input voltage $V_{in}$, duty-cycle $D$, load resistance $R_T$ and gate-driver voltage $V_{dr}$ and resistance $R_g$. Second, $f_s$, $V_{in}$, $D$ and $R_T$ are input *variables* of the switching power loss formula to be identified, whereas $V_{dr}$ and $R_g$ are *parameters* determining the values of the loss formula coefficients. All the quantities ($f_s$; $V_{in}$, $D$, $R_T$, $V_{dr}$, $R_g$) are easily determined for a given application.

III. THE GP-MOO MODELING APPROACH

Given an inverter PoM, a set of $m_1$ MOSFET gate driver voltage values $V_{dr,j1}$ has been considered, with $j_1 = 1,...,m_1$, and a set of $m_2$ MOSFET gate driver resistance values $R_{g,j2}$ has been analyzed, with $j_2 = 1,...,m_2$, thus resulting in the overall $m = m_1 \times m_2$ gate-driver conditions ($V_{dr,j}$, $R_{g,j}$), with $j = 1,...,m$. For each condition, $n$ combinations of switching frequency, input voltage, duty-cycle and total load resistance ($f_{s,i}$, $V_{in,i}$, $D_i$, $R_{T,i}$) have been considered, with $i = 1,...,n$. For each one of the $n \times m$ test conditions, a data vector has been created, including the test values ($f_{s,i}$; $V_{in,i}$, $D_i$, $R_{T,i}$) and the corresponding switching power loss value $y_{ij} = P_{sw}(f_{s,i}, V_{in,i}, D_i, R_{T,i}, V_{dr,j}, R_{g,j})$, simulated by means of the thermal-based analytical model described in Section II. In this paper, such $n \times m$ data vectors are the training data set $T$ used for switching power loss model identification. The main goal is to identify the behavioral model (11):

$$P_{sw,bhv} = \mathbf{F}\left[f_s, V_{in}, D, R_T, \mathbf{p}(V_{dr}, R_g)\right] \quad (11)$$

such that the value of the function $\mathbf{F}$ computed for each test condition of the training data set $T$ is as close as possible to the corresponding analytical model-based value $y_{ij}$, $\forall i \in \{1,...,n\}$ and $\forall j \in \{1,...,m\}$. In formula (11), $\mathbf{p}$ is a vector of numeric coefficients, given as a function of ($V_{dr}$, $R_g$). To discover this behavioral model, a GP algorithm has been considered [13]. A

GP is an evolutionary algorithm whose population is composed of "models". The population evolves, based on the standard genetic operations of selection, cross-over, mutation, elitism, etc. At the end of its evolution, the algorithm finds the models with the best-so-far fitness values. To construct the models in the population, the GP algorithm considers a set of functions (*non-terminal* set), and a set of constant coefficients and input variables (*terminal* set). Complexity factors $cf$ can be assigned to all the elements of these sets. Different combinations of $cf$ values have been investigated for this study. These $cf$ values have been eventually assumed for the terminal set: $cf = 0.6$ for multiplication of input variables, $cf = 1$ for all other operations between input variables and for constant coefficients. Instead, these $cf$ values have been considered for the non-terminal set: $cf = 1$ for sum and multiplication operators, $cf = 1.5$ for more complex functions, like logarithms, exponentials, powers, arctangents, hyperbolic tangents, etc.

Given the input variables ($f_s$; $V_{in}$, $D$, $R_T$) and the input parameters ($V_{dr}$, $R_g$) of the model, the structure of the behavioral power loss function **F** has to be the same for all the gate driver conditions ($V_{dr,j}$, $R_{g,j}$). Conversely, the model coefficients **p** are functions of ($V_{dr,j}$, $R_{g,j}$). To determine such coefficients for each gate-driver condition, a Non-Linear Least Squares (NLLS) algorithm, based on the Levenberg-Marquardt optimization method, has been applied to the respective $n$ experimental data vectors. Thus, the values of the coefficients **p** have been determined by minimizing the $\chi$-squared error between the analytical power loss $y_{ij}$ and the GP-predicted power loss $\mathbf{F}[f_{s,i}, V_{in,i}, D_i, R_{T,i}, \mathbf{p}(V_{dr,j}, R_{g,j})]$ for $i = 1,...,n$, as given in (12):

$$\chi_j^2 = \frac{1}{n}\sum_{i=1}^{n}\left\{\mathbf{F}\left[f_{s,i}, V_{in,i}, D_i, R_{T,i}, \mathbf{p}\left(V_{dr,j}, R_{g,j}\right)\right] - y_{ij}\right\}^2 \quad (12)$$

Then, interpolating functions $\mathbf{p}(V_{dr}, R_g)$ have been determined, as discussed hereafter. To select the best switching power loss model among all the discovered ones, the *accuracy* and the *complexity* of each GP-based model have been evaluated. For the model accuracy, the Root Mean Square Error (RMSE) between the analytical and the GP-predicted power losses has been evaluated over the entire training data set, as given in (13):

$$RMSE = \sqrt{\frac{1}{m}\sum_{j=1}^{m}\chi_j^2} \quad (13)$$

For the global complexity of each GP model, a term $F_{complexity}$ has been introduced and evaluated based on the complexity factors $cf$ of the elementary functions used within the model structure **F**. In particular, if a function is the argument of another function, then the complexity factors $cf$ of the two functions are multiplied. If two functions are multiplied or summed, then their complexity factors $cf$ are summed and multiplied by the complexity factor of a sum or a product, respectively. In the first case, a vertical development of the models (e.g., involved functions of functions) is prevented, especially for the functions with high $cf$ values. In the second case, a horizontal development of the models is avoided (e.g., models composed of many simple functions, multiplied or summed among them). Finally, an elitist Non-dominated Sorting Genetic Algorithm (NSGA-II [14]) has been adopted to discover the behavioral switching power loss model ensuring optimal trade-off between *RMSE* and $F_{complexity}$, selected as objective functions for minimization in this MOO problem.

## IV. RESULTS AND DISCUSSION

The discussion is herein referred to the VS-ETY020P120F part [15], a full-bridge inverter SiC MOSFET PoM from Vishay, whose nominal characteristics are $V_{ds} = 1200$V, $R_{ds,on} = 71$m$\Omega$, $I_d = 26$A. For other PoM parameters, the reader can refer to [15]. As a reference case study, the operating conditions of Table I have been considered. All the possible combinations of such values have been tested, resulting in a training data set composed of 1215 data vectors. The GP-MOO algorithm has been executed over 50 total independent runs, to verify the repeatability of the resultant behavioral models. The following metrics have been considered to classify each model:
- $N_{run}$ : number of GP runs during which the algorithm has discovered a certain model;
- $N_{gen}$ : average number of GP generations during which a model exists in the population;
- {$\mu_{err}$, $\sigma_{err}$, $err_{max}$}: mean value, standard deviation and maximum value of the distribution of the relative percent error provided by the GP model over the training data set.

In particular, for each test condition, a relative percent error of the GP model has been evaluated as given in (14):

$$err = \left\{\mathbf{F}\left[f_{s,i}, V_{in,i}, D_i, R_{T,i}, \mathbf{p}\left(V_{dr,j}, R_{g,j}\right)\right] - y_{ij}\right\}\frac{100}{y_{ij}} \quad (14)$$

The GP-MOO algorithm has eventually discovered 756 total models, characterized by different accuracy, complexity and metrics values. Table II summarizes the number of models characterized by different values of occurrences $N_{run}$ over 50 runs. Only the non-dominated Pareto-optimal solutions with $N_{run} \geq 6$ and $err_{max} \leq 80\%$ have been considered for comparisons. Fig. 3 shows such solutions with relevant {$N_{run}$, $N_{gen}$, $\mu_{err}$, $\sigma_{err}$, $err_{max}$} metrics values and model formulas. In particular, the model #2 given in (15) has optimal values of the metrics over the training data set, with $\mu_{err}=0.4\%$, $\sigma_{err}=3\%$, $err_{max}=12\%$:

TABLE I. OPERATING CONDITIONS FOR THE TRAINING DATA SET

| $f_s$ [kHz] | $V_{in}$ [V] | $D$ | $R_T$ [$\Omega$] | $V_{dr}$ [V] | $R_{gate}$ [$\Omega$] |
|---|---|---|---|---|---|
| 45; 75; 105 | 200; 300; 400 | 0.1; 0.3; 0.5; 0.7; 0.9 | 40; 70; 100 | 10; 15; 20 | 1; 3; 5 |

TABLE II. NUMBER OF GP-BASED MODELS WITH DIFFERENT $N_{run}$ VALUES

| $N_{run}$ | $1 \leq N_{run} < 2$ | $2 \leq N_{run} < 5$ | $5 \leq N_{run} < 10$ | $10 \leq N_{run} < 20$ | $20 \leq N_{run} < 51$ |
|---|---|---|---|---|---|
| # models | 613 | 101 | 24 | 8 | 10 |

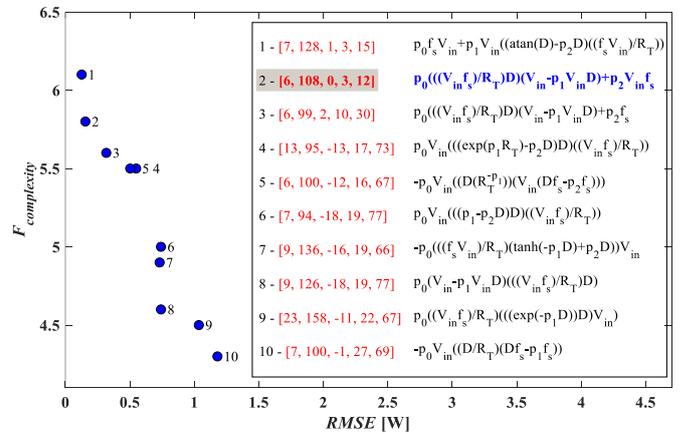

1 - [7, 128, 1, 3, 15]   $p_0 f_s V_{in} + p_1 V_{in}((\text{atan}(D) - p_2 D)((f_s V_{in})/R_T))$
2 - [6, 108, 0, 3, 12]   $p_0(((V_{in}f_s)/R_T)D)(V_{in} - p_1 V_{in}D) + p_2 V_{in}f_s$
3 - [6, 99, 2, 10, 30]   $p_0(((V_{in}f_s)/R_T)D)(V_{in} - p_1 V_{in}D) + p_2 f_s$
4 - [13, 95, -13, 17, 73]   $p_0 V_{in}(((\exp(p_1 R_T) - p_2 D)D)((V_{in}f_s)/R_T))$
5 - [6, 100, -12, 16, 67]   $-p_0 V_{in}((D(R_T^{-p_1}))(V_{in}(Df_s - p_2 f_s)))$
6 - [7, 94, -18, 19, 77]   $p_0 V_{in}((p_1 - p_2 D)D)((V_{in}f_s)/R_T))$
7 - [9, 136, -16, 19, 66]   $-p_0(((f_s V_{in})/R_T)(\tanh(-p_1 D) + p_2 D))V_{in}$
8 - [9, 126, -18, 19, 77]   $p_0(V_{in} - p_1 V_{in}D)(((V_{in}f_s)/R_T)D)$
9 - [23, 158, -11, 22, 67]   $p_0(V_{in}f_s)/R_T)((\exp(-p_1 D))DV_{in})$
10 - [7, 100, -1, 27, 69]   $-p_0 V_{in}((D/R_T)(Df_s - p_1 f_s))$

Fig. 3. Repeatable Pareto-optimal solutions

$$P_{sw,bhv} = p_0 f_s V_{in}^2 \cdot \frac{D(1-p_1 D)}{R_T} + p_2 f_s V_{in} \quad (15)$$

The coefficients $\{p_0, p_1, p_2\}$ are shown in Fig. 4 for different values of the gate-driver voltage $V_{dr}$ and for three gate-driver resistance values $R_g$. In particular, the coefficient $p_0$ decreases with $V_{dr}$ and increases with $R_g$, whereas the coefficient $p_2$ presents an opposite trend. The coefficient $p_1$ is nearly constant with $V_{dr}$ and $R_g$, being equal to about 1. Both $p_0$ and $p_2$ can be modeled as second-order polynomials of $V_{dr}$, according to (16):

$$p_k(V_{dr}) = a_0 V_{dr}^2 + a_1 V_{dr} + a_2 \quad \text{for } k = \{0, 2\} \quad (16)$$

A Linear Least Squares (LLS) algorithm has been used to determine the values of $\{a_0, a_1, a_2\}$ for each $R_g$ value. The resulting fitting curves of the coefficients $p_0$ and $p_2$ are shown in Fig. 4 (continuous lines). Linear and quadratic trends have been observed for the coefficients $\{a_0, a_1, a_2\}$ of $p_0$ and $p_2$, thus modeled as a function of $R_g$ according to (17):

$$a_x(R_g) = b_0 R_g^2 + b_1 R_g + b_2 \quad \text{for } x = \{0, 1, 2\} \quad (17)$$

where $b_0 = 0$ for the $\{a_0, a_1, a_2\}$ values of the coefficient $p_0$. Table III summarizes the $\{b_0, b_1, b_2\}$ values obtained by means of the LLS algorithm to fit the $\{a_0, a_1, a_2\}$ trends, for the coefficients $p_0$ and $p_2$.

Finally, the PoM switching power loss has been evaluated by means of (15)-(17), for all the combinations of the operating conditions given in Table I. Fig. 5 shows the relative percent errors between the results of the behavioral model $P_{sw,bhv}$ and the analytical model presented in Section II, for all the test conditions of the training data set $T$. Fig. 5 also shows the Probability Density Function (PDF) of the normal error distribution. The proposed behavioral model does provide reliable switching power loss estimation over the training data set, with percent errors within ±12%.

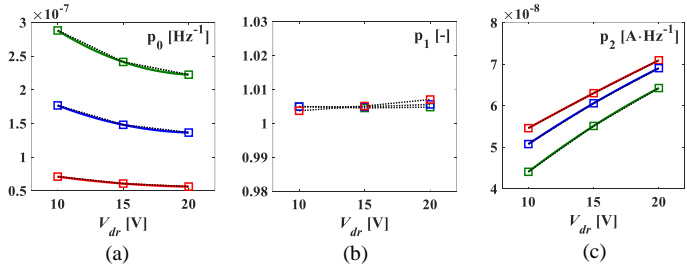

Fig. 4. Behavioral model coefficients and the relative fitting curves vs $V_{dr}$, for $R_g = 1\Omega$ (red), $R_g = 3\Omega$ (blue) and $R_g = 5\Omega$ (green).

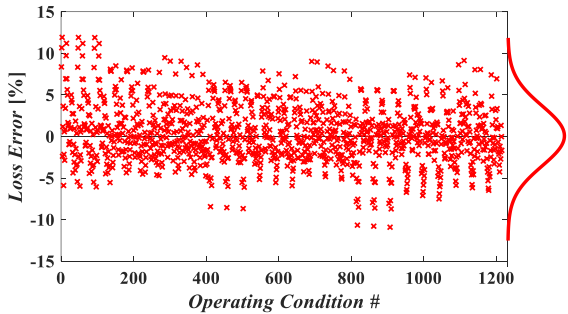

Fig. 5. Errors of the proposed behavioral model for training data set operating conditions and PDF of the normal error distribution.

TABLE III. FITTING CURVE COEFFICIENTS $\{b_0, b_1, b_2\}$

| Coefficients | | $b_0$ | $b_1$ | $b_2$ |
|---|---|---|---|---|
| $p_0$ | $a_0$ | 0 | 1.10e-10 | 9.62e-12 |
| | $a_1$ | 0 | -4.56e-09 | -5.02e-10 |
| | $a_2$ | 0 | 8.89e-08 | 1.99e-08 |
| $p_2$ | $a_0$ | 7.55e-13 | -1.16e-11 | 7.94e-13 |
| | $a_1$ | -2.33e-11 | 4.47e-10 | 1.51e-09 |
| | $a_2$ | -1.98e-10 | -3.80e-09 | 4.03e-08 |

CONCLUSIONS

This paper presents a new compact behavioral model for the evaluation of the switching power loss in full-bridge SiC MOSFET Power Modules (PoMs). A genetic programming algorithm and a multi-objective optimization approach have been jointly considered to identify a behavioral model as a function of the imposed operating conditions. The predictions of the proposed model have been validated over a wide range of operating conditions, by comparison to the switching power loss values calculated by using literature analytical models. As a prospective outcome, the proposed modeling procedure can be adopted by manufacturers to characterize their PoMs starting from experimental tests, thus providing power designers with simple, reliable and ready-to-use behavioral loss models.

ACKNOWLEDGEMENTS

The results here presented are developed in the framework of the 16ENG08 MICEV Project. The latter received funding from the EMPIR programme co-financed by the Participating States and from the European Union's Horizon 2020 research and innovation programme.